\documentclass[12pt,a4paper,oneside,reqno,notitlepage]{amsart}

\newtheorem{prop}{Proposition}
\newtheorem{theorem}{Theorem}

\usepackage{latexsym, amsbsy, amsmath, amsfonts, amssymb, amsthm, amscd} 
\begin{document}

\title{On certain questions of the free group automorphisms theory}

\author[Bardakov]{Valeriy G. Bardakov}
\address{Sobolev Institute of Mathematics, Novosibirsk 630090, Russia}
\email{bardakov@math.nsc.ru}

\author[Mikhailov]{Roman Mikhailov}
\address{Steklov Mathematical Institute, Gubkina 8, 119991 Moscow, Russia}
\email{rmikhailov@mail.ru}

\begin{abstract}
Certain subgroups of the groups $Aut(F_n)$ of automorphisms of a
free group $F_n$ are considered. Comparing Alexander polynomials
of two poly-free groups $Cb_4^+$ and $P_4$ we prove that these
groups are not isomorphic, despite the fact that they have a lot
of common properties. This answers the question of
Cohen-Pakianathan-Vershinin-Wu from \cite{CVW}. The questions of
linearity of subgroups of $Aut(F_n)$ are considered. As an
application of the properties of poison groups in the sense of
Formanek and Procesi, we show that the groups of the type
$Aut(G*\mathbb Z)$ for certain groups $G$ and the subgroup of
$IA$-automorphisms $IA(F_n)\subset Aut(F_n)$ are not linear for
$n\geq 3$. This generalizes the recent result of Pettet that
$IA(F_n)$ are not linear for $n\geq 5$.
\end{abstract}

\thanks{The research of the first author has been supported by the Scientific and
Research Council of Turkey (T${\rm \ddot{U}}$B${\rm \dot{I}}$TAK).
The research of the second author was partially supported by
Russian Foundation Fond of Fundamental Research, grant N
05-01-00993, and the Russian Presidential grant MK-1487.2005.1}
\maketitle

\section{Introduction}
\vspace{.5cm}
\subsection{The group of basis conjugating automorphisms} Let $F_n$ be a free group of rank $n\geq 2$ with a free generator
set $\{x_1, x_2,...,x_n\}$ and $Aut(F_n)$ the group of
automorphisms of $F_n$. Taking the quotient of $F_n$ by its
commutator subgroup $F_n',$ we get a natural homomorphism
$$
\xi : {\rm Aut}(F_n) \longrightarrow {\rm Aut}(F_n/F_n')={\rm
GL}_n(\mathbb Z),
$$
where ${\rm GL}_n(\mathbb Z)$ is the general linear group over the
ring of integers. The kernel of this homomorphism consists of
automorphisms acting trivially modulo the commutator subgroup
$F'_n$. It is called {\it the group of $IA$-automorphisms} and
denoted by ${\rm IA}(F_n)$ (see \cite[Chapter 1, \S~4]{LS}). The
group ${\rm IA}(F_2)$ is isomorphic to the group of inner
automorphisms ${\rm Inn}(F_2)$, which is isomorphic to the free
group $F_2$.

J. Nilsen (for the case $n\leq 3$) and W. Magnus (for all $n$)
have shown that (see \cite[Chapter 1, \S~4]{LS}) the group ${\rm
IA}(F_n)$ is generated by the following automorphisms
$$
\varepsilon_{ijk} : \left\{
\begin{array}{ll}
x_{i} \longmapsto x_{i}[x_j, x_k] & \mbox{if }~~ k\neq i, j, \\
x_{l} \longmapsto x_{l} & \mbox{if }~~ l\neq i,
\end{array} \right. ~~~~~
\varepsilon_{ij} : \left\{
\begin{array}{ll}
x_{i} \longmapsto x_{j}^{-1}x_ix_j & \mbox{if }~~ i\neq j, \\
x_{l} \longmapsto x_{l} & \mbox{if }~~ l\neq i.
\end{array} \right.
$$

The subgroup of the group ${\rm IA}(F_n)$ is generated by the
automorphisms $\varepsilon_{ij}$, $1\leq i\neq j\leq n$ is called
{\it the group of basis conjugating automorphisms}. We will denote
this group by $Cb_n$. The group $Cb_n$ is a subgroup of {\it the
group of conjugating automorphisms} $C_n$. (Recall that any
automorphism of $C_n$ sends a generator $x_i$ to an element of the
type $f_i^{-1}x_{\pi (i)}f_i$, where $f_i \in F_n$, and $\pi$ is a
permutation from the symmetric group $S_n$.) Clearly, if $\pi$ is
the identity permutation, then the described element lies in the
group $Cb_n$.

It is shown by J. McCool \cite{Mc} that the group of basis
conjugating automorphisms $Cb_n$ is generated by the automorphisms
$\varepsilon_{ij}$, $1\leq i\neq j \leq n$, defined above has the
following relations (we denote different indexes by different
symbols):
$$
\begin{array}{lr}
\varepsilon_{ij}\varepsilon_{kl}=\varepsilon_{kl}\varepsilon_{ij},
&  \\
& \\
\varepsilon_{ij}\varepsilon_{kj}=\varepsilon_{kj}\varepsilon_{ij},
&   \\
& \\
(\varepsilon_{ij}\varepsilon_{kj})\varepsilon_{ik}=\varepsilon_{ik}(\varepsilon_{ij}
\varepsilon_{kj}).
& \\
\end{array}
$$

It is shown in \cite{Bar} that the group of conjugating
automorphisms $Cb_n$, $n\geq 2$ decomposes into a semi-direct
product:
\begin{equation}\label{razl1}
Cb_n=D_{n-1}\leftthreetimes (D_{n-2}\leftthreetimes (\ldots
\leftthreetimes (D_2\leftthreetimes D_1))\ldots ),
\end{equation}
 where the group $D_i,$ $i=1,2,\ldots,n-1,$ is generated by the elements
 $\varepsilon_{i+1,1},\varepsilon_{i+1,2},
\ldots,\varepsilon_{i+1,i},$ $\varepsilon_{1,i+1},$
$\varepsilon_{2,i+1},$ $\ldots,\varepsilon_{i,i+1}$.
 The elements $\varepsilon_{i+1,1},\varepsilon_{i+1,2},
\ldots,\varepsilon_{i+1,i}$ generate a free group of rank $i$, but
the elements
$\varepsilon_{1,i+1},\varepsilon_{2,i+1},\ldots,\varepsilon_{i,i+1}$
generate a free abelian group of rank $i$.

\subsection{The braid group as a subgroup of  $Aut(F_n)$}
The braid group $B_n$, $n\geq 2$, on $n$ strands is defined by the
generators $\sigma_1,\sigma_2,\ldots,\sigma_{n-1}$ with defining
relators
$$
\sigma_i\sigma_{i+1}\sigma_i=\sigma_{i+1}\sigma_i\sigma_{i+1}~~~\mbox{if
}~~ i=1,2,\ldots,n-2, 
$$
$$
\sigma_i \sigma_j = \sigma_j \sigma_i~~~\mbox{if }~~ |i-j|\geq 2.
$$
There exists a natural homomorphism from the braid group $B_n$
onto the permutation group $S_n$, sending the generator $\sigma_i$
to the transposition  $(i,i+1)$, $i=1,2,\ldots,n-1$. The kernel of
this homomorphism is called {\it the pure braid group} and denoted
by $P_n$. The group $P_n$ is generated by the elements $a_{ij}$,
$1\leq i < j\leq n$ which can be expressed in terms of generators
of $B_n$ as follows:
$$
a_{i,i+1}=\sigma_i^2,
$$
$$
a_{ij}=\sigma_{j-1}\sigma_{j-2}\ldots\sigma_{i+1}\sigma_i^2\sigma_{i+1}^{-1}\ldots
\sigma_{j-2}^{-1}\sigma_{j-1}^{-1},~~~i+1< j \leq n.
$$
The pure braid group $P_n$ is a semi-direct product of the normal
subgroup $U_n$, which is the free group with free generators
$a_{1n}, a_{2n},\ldots,a_{n-1,n}$ with the group $P_{n-1}$.
Analogously, $P_{n-1}$ is a semi-direct product of the free
subgroup $U_{n-1}$ with the free generators $a_{1,n-1},
a_{2,n-1},\ldots,a_{n-2,n-1}$ with the subgroup $P_{n-2}$ etc.
Hence the group $P_n$ has the following decomposition
\begin{equation}\label{razl2}
P_n=U_n\leftthreetimes (U_{n-1}\leftthreetimes (\ldots
\leftthreetimes (U_3\leftthreetimes U_2))\ldots),~~~U_i\simeq
F_{i-1}, i=2,3,\ldots,n.
\end{equation}
The pure braid group  $P_n$ is defined by relations (for $\nu=\pm
1$)
$$
\begin{array}{ll}
a_{ik}^{-\nu }a_{kj}a_{ik}^{\nu }=\left( a_{ij}a_{kj}\right) ^{\nu
}a_{kj}
\left( a_{ij}a_{kj}\right) ^{-\nu }, & \\
& \\
a_{km}^{-\nu }a_{kj}a_{km}^{\nu }=\left( a_{kj}a_{mj}\right) ^{\nu
}a_{kj}
\left( a_{kj}a_{mj}\right) ^{-\nu }, &~~~m < j,\\
& \\
a_{im}^{-\nu }a_{kj}a_{im}^{\nu }=\left[ a_{ij}^{-\nu },
a_{mj}^{-\nu }\right] ^{\nu }a_{kj} \left[ a_{ij}^{-\nu },
a_{mj}^{-\nu }\right] ^{-\nu }, & ~~~i < k <
m,\\
& \\
a_{im}^{-\nu }a_{kj}a_{im}^{\nu }=a_{kj}, & ~~~k < i;~~~m < j~~
\mbox{or}~~ m < k,
\end{array}
$$
where $[a, b]=a^{-1}b^{-1}ab$ is the commutator of elements $a$
and $b$.

The braid group $B_n$ can be embedded into the automorphism group
 ${\rm Aut}(F_n)$. In this embedding the generator $\sigma_i$,
$i=1,2,\ldots,n-1$, defines the following automorphism
$$
\sigma_{i} : \left\{
\begin{array}{ll}
x_{i} \longmapsto x_{i}x_{i+1}x_i^{-1}, &  \\ x_{i+1} \longmapsto
x_{i}, & \\ x_{l} \longmapsto x_{l} & \mbox{if }~~ l\neq i,i+1.
\end{array} \right.
$$
The generator $a_{rs}$ of the pure braid group $P_n$ defines the
following automorphism
$$
a_{rs} : \left\{
\begin{array}{ll}
x_{i} \longmapsto x_{i} & \mbox{if }~~ s < i~~ \mbox{or}~~ i < r, \\
& \\
x_{r} \longmapsto x_{r}x_{s}x_rx_s^{-1}x_r^{-1}, &  \\
& \\
x_{i} \longmapsto [x_{r}^{-1}, x_s^{-1}]x_i[x_{r}^{-1}, x_s^{-1}]^{-1} & \mbox{if }~~ r < i <s,\\
& \\
x_{s} \longmapsto x_{r}x_sx_r^{-1}. &
\end{array} \right.
$$

As was shown by E. Artin, the automorphism $\beta $ from ${\rm
Aut}(F_n)$ belongs to the braid group $B_n$ if and only if $\beta
$ satisfies the following two conditions:
\begin{align*} & 1)\
\beta(x_i)=a_i^{-1}x_{\pi(i)}a_i,~~~1\leq i\leq n,\\
& 2)\ \beta(x_1x_2 \ldots x_n)=x_1x_2 \ldots x_n,
\end{align*}
where  $\pi $ is a permutation from $S_n$, and $a_i\in F_n$.

The set of automorphisms from $C_n$, which act trivially on the
product $x_1x_2...x_n$ is exactly  the braid group $B_n$. It was
shown by A. Savushkina \cite{S}, that the braid group $B_n$
intersects with the subgroup $Cb_n$ by the pure braid group $P_n$.
Thus the group $P_n$ is a subgroup of the basis conjugating
automorphism group $Cb_n$ for all $n\geq 2$. Furthermore, the
concordance of decompositions  (\ref{razl1}) and (\ref{razl2}) of
groups $Cb_n$ and $P_n$ respectively, takes a place: there are
corresponding embeddings  $U_{i+1}\leq D_i,$ $i=1,2,\ldots,n-1.$

\subsection{The group $Cb_n^+$}
Denote by $Cb_n^+$ the subgroup of the group $Cb_n,$ $n \geq 2$,
generated by the elements $\varepsilon_{ij},\ i>j$. Groups
$Cb_n^+$ are poly-free groups: it is shown in \cite{CCP} that
there exists the splitting exact sequence:
$$
1\to F_{n-1}\to Cb_n^+\to Cb_{n-1}^+\to 1,
$$
where $F_{n-1}$ is the free subgroup in $Cb_n^+$, generated by the
elements
$$
\varepsilon_{n1}, \varepsilon_{n2}, \ldots, \varepsilon_{n,n-1}.
$$

Description of the Lie algebra constructed from lower central
filtration of the basis conjugating automorphism group, is a
non-trivial problem; the methods used in the description of Lie
algebras for pure braid groups, do not work in this case. However,
in the case of the group $Cb_n^+$, the situation is different. Lie
algebras and cohomology rings of groups $Cb_n^+$ were described in
\cite{CVW}. It was also shown in \cite{CVW} that the group
$Cb_n^+$ is isomorphic to the pure braid group $P_n$ for $n=2,3$.
Groups $Cb_n^+$ and $P_n$ are quite similar. For instance, the
groups $Cb_n^+$ and $P_n$ are {\it stably isomorphic}, namely,
suspensions over their classifying spaces $\Sigma K(Cb_n^+,1)$ and
$\Sigma K(P_n,1)$ are homotopically equivalent for all  $n\geq 1$
(see \cite{CVW}). The conjecture that $Cb_n^+$ and $P_n$ are
isomorphic for all $n\geq 1$ comes naturally (question 1,
\cite{CVW}). One of the main results of this paper is the proof
that the groups $Cb_4^+$ and $P_4$ are {\it non-isomorphic}
(Theorem \ref{the_on_iso}).
\subsection{The questions of linearity}
As it was shown above, for any $n\geq 2$, there is the following
chain of subgroups of the automorphism group $Aut(F_n)$:
$$
P_n\subset Cb_n\subset IA_n(F_n)\subset Aut(F_n).
$$
It was shown in \cite{Bigelow}, \cite{Krammer} that the braid
groups $B_n$ are linear. Hence the pure braid groups $P_n$ are
also are linear. The situation with groups $Aut(F_n)$ is as
follows. The group ${\rm Aut}(F_n)$ is not linear for all $n\geq
3$ (see \cite{FP}), however the group ${\rm Aut}(F_2)$ is linear
(such a representation was constructed, for example, in
\cite{Bar2}). A certain non-linear subgroup (called {\it poison
group}) of ${\rm Aut}(F_n)$ was constructed in \cite{FP}. The
following question rises naturally (see question 15.17 \cite{KT}):
are the groups of IA-automorphisms ${\rm IA}(F_n)$ and the groups
of basis conjugating automorphisms  $Cb_n$ linear for $n\geq 2$?
Recently A. Pettet showed that the group ${\rm IA}(F_n)$ is not
linear for $n \geq 5$. In this paper we show that the groups ${\rm
IA}(F_n)$ are not linear for all $n \geq 3$ (Theorem
\ref{cor_on_lineal}). Therefore, the complete answer to the
question 15.17 from \cite{KT} is given.

\section{The group $Cb_4^+$ is not isomorphic to the pure braid group $P_4$}
\vspace{.5cm} \subsection{Fox differential calculus} Recall the
definition and main properties of Fox derivatives \cite[chapter
3]{B}, \cite[chapter 7]{CF}.

Let $F_n$ be a free group of rank $n$ with free generators $x_1,
x_2, \ldots, x_n$. Let $\varphi $ be an endomorphism of the group
$F_n$. Denote by  $F_n^{\varphi}$ the image of $F_n$ under the
endomorphism $\varphi$. Let $\mathbb{Z}$ be the ring of integers,
$\mathbb{Z}G$ the integral group ring of  $G$.

For every $j=1, 2, \ldots, n$ define a map
$$
\frac{\partial }{\partial x_j} : \mathbb{Z}F_n\longrightarrow
\mathbb{Z}F_n
$$
by setting
$$
1)~~~ \frac{\partial x_i}{\partial x_j}= \left\{
\begin{array}{ll}
1 & \mbox{for}\ i = j , \\
0 & \mbox{otherwise},
\end{array} \right.
$$
$$
2)~~~ \frac{\partial x_i^{-1}}{\partial x_j}= \left\{
\begin{array}{ll}
-x_i^{-1} & \mbox{for}\ i = j , \\
0 & \mbox{otherwise},
\end{array} \right.
$$
$$
3)~~~ \frac{\partial (w v)}{\partial x_j} = \frac{\partial
w}{\partial x_j} (v)^{\tau } + w \frac{\partial v}{\partial
x_j},~~~w, v \in \mathbb{Z}F_n,
$$
where $\tau : \mathbb{Z}F_n \longrightarrow \mathbb{Z}$ is the
trivialization operation, sending all elements of $F_n$ to the
identity,
$$
4)~~~ \frac{\partial }{\partial x_j} \left( \sum a_g g\right) = \sum
a_g \frac{\partial g}{\partial x_j},~~~g\in F_n,~~~a_g \in
\mathbb{Z}.
$$

Denote by $\Delta_n$ the augmentation ideal of the ring
$\mathbb{Z}F_n$, i.e. the kernel of the homomorphism $\tau$. It is
easy to see that for every $v \in \mathbb{Z}F_n$, the element
$v-v^{\tau}$ belongs to $\Delta_n$. Furthermore, there is the
following formula
\begin{equation}\label{fc}
v-v^{\tau } = \sum\limits_{j=1}^n \frac{\partial v}{\partial x_j}
(x_j - 1),
\end{equation}
which is called {\it the main  formula of Fox calculus}. Formula
(\ref{fc}) implies that the elements  $\{x_1 - 1, x_2 - 1, \ldots,
x_n - 1\}$ determine the basis of the augmentation ideal
$\Delta_n$.

\subsection{Groups $Cb_4^+$ and $P_4$} As it was noted in
the introduction, the poly-free groups $Cb_n^+$ and $P_n$ have
similar properties. Clearly, one has
\begin{align*}
& Cb_2^+\simeq P_2\simeq \mathbb Z;\\
& Cb_3^+\simeq P_3\simeq F_2\rtimes \mathbb Z.
\end{align*}
However, the next result shows that the groups $Cb_4^+$ and $P_4$
being semidirect products of $F_3$ with $F_2\rtimes \mathbb Z$,
are not isomorphic.
\begin{theorem}\label{the_on_iso}
Groups $Cb_4^+$ and $P_4$ are not isomorphic.
\end{theorem}
\begin{proof}
It is well-known  (see \cite{Ne}) that the group $P_4$ decomposes
as a direct product of its center generated by a single element
$a_{12} a_{13} a_{23} a_{14} a_{24} a_{34}$ and a subgroup $H=U_4
\leftthreetimes U_3,$ $U_4 = \langle a_{14}, a_{24},
a_{34}\rangle,$ $U_3 = \langle a_{13}, a_{23} \rangle$. It is
easy to check that the center of the group $Cb_4^+$
 is the infinite cyclic group generated by the element $\varepsilon_{21} \varepsilon_{31}
\varepsilon_{41}$ and the whole group $Cb_4^+$ is a direct product
of its center and a subgroup $G = D_3^+ \leftthreetimes D_2^+$,
where $D_3^+ = \langle \varepsilon_{41}, \varepsilon_{42},
\varepsilon_{43} \rangle, $ $D_2^+ = \langle \varepsilon_{31},
\varepsilon_{32} \rangle$. Since a center is a characteristic
subgroup for every group, it is enough to show that the group $G$
is not isomorphic to the group $H$. We will compare Alexander
polynomials for the groups $G$ and $H$.

The presentation of the group $Cb_4$ implies that the group $G$ is
defined by the following set of relations:
\begin{align*}
& \varepsilon_{31}^{-1} \varepsilon_{41} \varepsilon_{31} =
\varepsilon_{41},~~~ \varepsilon_{31}^{-1} \varepsilon_{42}
\varepsilon_{31} = \varepsilon_{42},~~~ \varepsilon_{31}^{-1}
\varepsilon_{43} \varepsilon_{31} = \varepsilon_{41}
\varepsilon_{43} \varepsilon_{41}^{-1},\\
& \varepsilon_{32}^{-1} \varepsilon_{41} \varepsilon_{32} =
\varepsilon_{41},~~~ \varepsilon_{32}^{-1} \varepsilon_{42}
\varepsilon_{32} = \varepsilon_{42},~~~ \varepsilon_{32}^{-1}
\varepsilon_{43} \varepsilon_{32} = \varepsilon_{42}
\varepsilon_{43} \varepsilon_{42}^{-1}.
\end{align*}
We present these relations in the following form:
\begin{align*} &
r_{11} = \varepsilon_{41}^{-1} \varepsilon_{31}^{-1}
\varepsilon_{41} \varepsilon_{31},~~~ r_{21} =
\varepsilon_{42}^{-1} \varepsilon_{31}^{-1} \varepsilon_{42}
\varepsilon_{31},~~~ r_{31} = \varepsilon_{43}^{-1}
\varepsilon_{41}^{-1} \varepsilon_{31}^{-1}
\varepsilon_{43} \varepsilon_{31} \varepsilon_{41},\\
& r_{12} = \varepsilon_{41}^{-1} \varepsilon_{32}^{-1}
\varepsilon_{41} \varepsilon_{32},~~~ r_{22} = \varepsilon_{42}^{-1}
\varepsilon_{32}^{-1} \varepsilon_{42} \varepsilon_{32},~~~ r_{32} =
\varepsilon_{43}^{-1} \varepsilon_{42}^{-1} \varepsilon_{32}^{-1}
\varepsilon_{43} \varepsilon_{32} \varepsilon_{42}.
\end{align*}

We define the homomorphism $\varphi: G \to \langle t \rangle$, by
setting
$$\varphi(\varepsilon_{31}) = \varphi(\varepsilon_{32}) = \varphi(\varepsilon_{41}) =
\varphi(\varepsilon_{42}) = \varphi(\varepsilon_{43}) = t
$$
and extend it to the homomorphism of integral group rings by
linearity
$$
\varphi : \mathbb{Z} G \longrightarrow \mathbb{Z} \langle t \rangle.
$$

Let us find Fox derivatives of the relators $r_{ij},$ with respect
to the symbols $\varepsilon_{31},$ $\varepsilon_{32},$
$\varepsilon_{41},$ $\varepsilon_{42},$ $\varepsilon_{43}$ and
their $\varphi$-images. We have the following
\begin{align*}
& \left( \frac{\partial r_{11}}{\partial \varepsilon_{31}}
\right)^{\varphi} =t^{-2}(t - 1), \left( \frac{\partial
r_{11}}{\partial \varepsilon_{41}} \right)^{\varphi} =t^{-2}(1 -
t), \left( \frac{\partial r_{12}}{\partial \varepsilon_{32}}
\right)^{\varphi} =t^{-2}(t - 1),\\
& \left( \frac{\partial r_{12}}{\partial \varepsilon_{41}}
\right)^{\varphi} =t^{-2}(1 - t), \left( \frac{\partial
r_{21}}{\partial \varepsilon_{31}} \right)^{\varphi} =t^{-2}(t -
1), \left( \frac{\partial r_{21}}{\partial \varepsilon_{42}}
\right)^{\varphi} =t^{-2}(1 - t),\\
& \left( \frac{\partial r_{22}}{\partial \varepsilon_{32}}
\right)^{\varphi} =t^{-2}(t - 1), \left( \frac{\partial
r_{22}}{\partial \varepsilon_{42}} \right)^{\varphi} =t^{-2}(1 -
t), \left( \frac{\partial r_{31}}{\partial \varepsilon_{31}}
\right)^{\varphi} =t^{-1}(t - 1),\\
& \left( \frac{\partial r_{31}}{\partial \varepsilon_{41}}
\right)^{\varphi} =t^{-2}(t - 1), \left( \frac{\partial
r_{31}}{\partial \varepsilon_{43}} \right)^{\varphi} =t^{-3}(1 -
t^2), \left( \frac{\partial r_{32}}{\partial \varepsilon_{32}}
\right)^{\varphi} =t^{-1}(t - 1),\\
& \left( \frac{\partial r_{32}}{\partial \varepsilon_{42}}
\right)^{\varphi} =t^{-2}(t - 1), \left( \frac{\partial
r_{32}}{\partial \varepsilon_{43}} \right)^{\varphi} =t^{-3}(1 -
t^2).
\end{align*}
Derivatives with respect to other generators are zero.

With the obtained values of derivatives we form the Alexander
matrix. After elementary transformations of rows and columns and
deleting zero rows and columns, we get the diagonal matrix
$$
(1 - t)\ {\rm diag(-t^{-2}, -t^{-2}, -t^{-2}, -t^{-2} (1+t))}.
$$
Since the Alexander polynomial is defined up to the multiplication
with a unit of the ring $\mathbb{Z} \langle t \rangle$, we get
$$
\Delta_G (t) = (1 - t)^4 (1 + t).
$$

\medskip

We now consider the group $H$. Defining relations of the group
$P_4$ imply that the group $H$ can be defined by the following
relations
\begin{align*}
& a_{13} a_{14} a_{13}^{-1} = a_{34}^{-1} a_{14} a_{34},~~~
a_{13}^{-1} a_{24} a_{13} = [a_{14}^{-1}, a_{34}^{-1}] a_{24}
[a_{34}^{-1}, a_{14}^{-1}],\\
& a_{23}^{-1} a_{14} a_{23} =  a_{14},  a_{23} a_{24} a_{23}^{-1} =
a_{34}^{-1} a_{24} a_{34},\\
& a_{13}^{-1} a_{34} a_{13} = a_{14} a_{34} a_{14}^{-1},~~~
a_{23}^{-1} a_{34} a_{23} = a_{24} a_{34} a_{24}^{-1}.
\end{align*}

These relations can be presented in the following form
\begin{align*} & q_{11} = a_{13} a_{14} a_{13}^{-1} a_{34}^{-1}
a_{14}^{-1} a_{34},~~~ q_{21} = a_{13}^{-1} a_{24} a_{13} a_{14}
a_{34} a_{14}^{-1} a_{34}^{-1}
a_{24}^{-1} a_{34} a_{14} a_{34}^{-1} a_{14}^{-1},\\
& q_{12} =
a_{14}^{-1} a_{23}^{-1} a_{14} a_{23},~~~ q_{22} = a_{23}
a_{24} a_{23}^{-1} a_{34}^{-1} a_{24}^{-1} a_{34},\\
& q_{31} = a_{13}^{-1} a_{34} a_{13} a_{14} a_{34}^{-1}
a_{14}^{-1},~~~ q_{32} = a_{23}^{-1} a_{34} a_{23} a_{24}
a_{34}^{-1} a_{24}^{-1}.
\end{align*}

As above, we define the homomorphism $\varphi$ from the group $H$
to the infinite cyclic group $\langle t \rangle$, by setting
$$\varphi(a_{13}) = \varphi(a_{23}) = \varphi(a_{14}) =
\varphi(a_{24}) = \varphi(a_{34}) = t,
$$
and extend it by linearity to the homomorphism of group rings
$$
\varphi : \mathbb{Z} H \longrightarrow \mathbb{Z} \langle t \rangle.
$$

We compute the Fox derivatives of the relators $q_{ij},$ written
above, with respect to the variables  $a_{13},$ $a_{23},$
$a_{14},$ $a_{24},$ $a_{34}$ and find their images under the
homomorphism $\varphi$. We get the following identities (the
derivatives of all other relators can be written analogously)
\begin{align*}
& \left( \frac{\partial q_{11}}{\partial a_{13}} \right)^{\varphi}
= 1 - t, \left( \frac{\partial q_{11}}{\partial a_{14}}
\right)^{\varphi} = t^{-1} (t^2 - 1), \left( \frac{\partial
q_{11}}{\partial a_{34}} \right)^{\varphi} = t^{-1} (1 - t),\\
& \left( \frac{\partial q_{21}}{\partial a_{13}} \right)^{\varphi}
= t^{-1} (t - 1), \left( \frac{\partial q_{21}}{\partial a_{14}}
\right)^{\varphi} = -(1 - t)^2, \left( \frac{\partial
q_{21}}{\partial a_{24}} \right)^{\varphi} = t^{-1} (1 - t),\\
& \left( \frac{\partial q_{21}}{\partial a_{34}} \right)^{\varphi}
= (t - 1)^2, \left( \frac{\partial q_{32}}{\partial a_{23}}
\right)^{\varphi} = t^{-1} (t - 1), \left( \frac{\partial
q_{32}}{\partial a_{24}} \right)^{\varphi} = t - 1,\\ & \left(
\frac{\partial q_{32}}{\partial a_{34}} \right)^{\varphi} = t^{-1}
(1 - t^2).
\end{align*}

We form the matrix from the values of the calculated Fox
derivatives. After elementary transformations over rows and
columns, we get the diagonal matrix
$$ (1 - t)\ {\rm diag(1, -t^{-2}, 1, t^{-2} +
t^{-1} + 1)}.
$$
The Alexander polynomial of the group $H$ is equal to $$ \Delta_H
(t) = (1 - t)^4 (t^2 + t + 1).
$$
Thus the Alexander polynomials of the groups $G$ and $H$ are
different, and hence, the groups $G$ and $H$ are not isomorphic
(see, for example \cite{CF} for the proof that Alexander
invariants, in particular the Alexander polynomial, are invariants
of a group and do not depend on a given group presentation).
\end{proof}

\section{Isomorphism problem for a certain class of Lie algebras}
\vspace{.5cm}
\subsection{A class of  Lie algebras} The last section gives a motivation for
finding methods for proving that two given subgroup of
automorphism groups are not isomorphic. In this section we will
consider a similar question.

Consider a homomorphism
$$
\phi: F_n \to IA(F_k),
$$
where $n, k\geq 2$. We can form a natural semi-direct product
$$
G_\phi:=F_k\rtimes F_n
$$
and ask whether $G_{\phi_1}$ is isomorphic to $G_{\phi_2}$ for two
different homomorphisms $\phi_1,\phi_2: F_n\to IA(F_k)$.

For a given group $G$, denote by $L(G)$ the Lie algebra
constructed from the lower central series filtration:
\begin{align*} & L(G):=\bigoplus_{i\geq
1}\gamma_i(G)/\gamma_{i+1}(G),\\ & \gamma_1(G)=G,\
\gamma_{i+1}(G)=[\gamma_i(G),G]. \end{align*}

Here we will consider the class of Lie algebras $L(G_\phi)$ for
different homomorphisms $\phi$ and prove that some Lie algebras
from this class are not isomorphic.

Since $F_n$ acts trivially on $H_1(F_k)$, there is the following
exact splitting sequence of Lie algebras:
$$
0\to L(F_k)\to L(G_\phi)\to L(F_n)\to 0
$$
by \cite{FR}.
\subsection{Scheuneman's invariants}
Let $k$ be a field and $L$ a finitely-generated Lie algebra of
nilpotence class two. Let us choice a basis of $L$: $\{x_1,\dots,
x_n, y_1,\dots, y_r\}$ such that $\{y_1,\dots, y_r\}$ is a basis
of $[L,L]$. Let $U(L)$ be a universal enveloping algebra of $L$.
Then $U(L)$ can be identified with a polynomial algebra over $k$
with non-commutative variables $x_1,\dots, x_n,y_1,\dots, y_r$ and
relations \begin{align*} & x_iy_j=y_jx_i,\ \text{for all}\ i,j;\\
& y_iy_j=y_jy_i,\ \text{for all}\ i,j;\\
& x_ix_j-x_jx_i=\text{linear combination of}\ y_k\text{'s};
\end{align*}
Consider an alternative sum
\begin{equation}\label{poly}
I(x_1,\dots, x_n)=\sum_{\sigma=(i_1,\dots, i_n)\in
S_n}(-1)^{|\sigma|}x_{i_1}\dots x_{i_n}, \end{equation} where
$|\sigma|$ denotes the sign of the permutation $\sigma$.

Two polynomials $f(z_1,\dots, z_m)$ and $g(z_1,\dots, z_m)$ are
{\it $k$-equivalent} if $$ f(z_1,\dots, z_m)=ag(z_1',\dots, z_m'),
$$
where $a\in k, a\neq 0$,
$$
z_j'=\sum b_{ij}z_j,\ b_{ij}\in k,\ det(b_{ij})\neq 0.
$$
The following theorem is due to Scheuneman:
\begin{theorem}\label{sch}\cite{Sch}
The $k$-equivalence class of $I(x_1,\dots,x_n)$ is an invariant of
$k$-isomorphism of the Lie algebra $L$.
\end{theorem}
For $k$-equivalent polynomials $f(z_1,\dots,z_m)$ and
$g(z_1,\dots, z_m)$ and $z_i'=\sum b_{ij}z_j,\ b_{ij}\in k,$
 $det(b_{ij})\neq 0,$ such that $f(z_1,\dots, z_m)=ag(z_1',\dots,
z_m'),$
$$
Hes(f)(z_1,\dots, z_m)=a^m (det(b_{ij}))^2Hes(g)(z_i',\dots,
z_m'),
$$
where $Hes(f)(z_1,\dots,z_m)$ is the Hessian
$det(\frac{\partial^2f}{\partial z_i\partial z_j})$ (see, for
example, Lemma 8 \cite{Sch}).

\subsection{Non-isomorphic Lie algebras}
Consider the following homomorphism $\phi_\alpha$ of a free group
$F_3$ with basis $\{u_1,u_2,u_3\}$ to the automorphism group of a
free group $F_3$ with basis $\{t_1,t_2,t_3\}$:
\begin{align*}
& \phi_{\alpha}: u_1\mapsto \left\{
\begin{array}{ll}
t_1\longmapsto t_1\alpha_1\\
t_i \longmapsto t_i & \mbox{for}\ t=2,3,
\end{array} \right.\\
& \phi_{\alpha}: u_2\mapsto \left\{
\begin{array}{ll}
t_2\longmapsto t_2\alpha_2\\
t_i \longmapsto t_i & \mbox{for}\ t=1,3,
\end{array} \right.\\
& \phi_{\alpha}: u_3\mapsto \left\{
\begin{array}{ll}
t_3\longmapsto t_3\alpha_3\\
t_i \longmapsto t_i & \mbox{for}\ t=1,2,
\end{array} \right.
\end{align*}
where $\alpha=(\alpha_1,\alpha_2,\alpha_3),$
\begin{align*}& \alpha_1\in \{[t_1,t_2],[t_2,t_3],[t_1,t_3]\},\\
& \alpha_2\in \{[t_2,t_1],[t_2,t_3],[t_1,t_3]\},\\ & \alpha_3\in
\{[t_1,t_2],[t_3,t_2],[t_3,t_1]\}.\end{align*}

Lie algebras $L(G_\alpha)$ which correspond to the groups
$G_\alpha$ can be described as Lie algebras with generators
$t_1,t_2,t_3,u_1,u_2,u_3$ and defining relations:
\begin{align*}
& [t_1,u_1]=\alpha_1,\ [t_2,u_2]=\alpha_2,\ [t_3,u_3]=\alpha_3,\\
& [t_2,u_1]=[t_3,u_1]=[t_1,u_2]=[t_3,u_2]=[t_1,u_3]=[t_2,u_3]=0.
\end{align*}

Consider the quotient
$$
L_\alpha:=L(G_{\phi_\alpha})/[L(G_{\phi_\alpha}),L(G_{\phi_\alpha}),L(G_{\phi_\alpha})].
$$
Clearly, $L_\alpha$ is 12-dimensional algebra with basis $$
\{t_1,t_2,t_3,u_1,u_2,u_3,y_1,y_2,y_3,y_4,y_5,y_6\},
$$
where $y_1=[t_1,t_2], y_2=[t_1,t_3], y_3=[t_2,t_3], y_4=[u_1,u_2],
y_5=[u_1,u_3], y_6=[u_2,u_3]$. Consider the polynomial
$I(t_1,t_2,t_3,u_1,u_2,u_3)$. The direct computation gives the
following expression
\begin{multline*}
I(t_1,t_2,t_3,u_1,u_2,u_3) =\\
[t_1,t_2][t_3,u_3][u_1,u_2]-[t_1,t_3][t_2,u_2][u_1,u_3]+[t_1,u_1][t_2,t_3][u_2,u_3]
-[t_1,u_1][t_2,u_2][t_3,u_3]\\
=y_1y_4\alpha_3-y_2y_5\alpha_2+y_3y_6\alpha_1-\alpha_1\alpha_2\alpha_3=P_{(\alpha_1,\alpha_2,\alpha_3)}(y_1,\dots,
y_6).
\end{multline*}
By Theorem \ref{sch}, we have
$$
\left\{ k-\text{equivalence classes of cubic forms}\
P_{(\alpha_1,\alpha_2,\alpha_3)}\}\right. \subseteq
\left\{k-\text{isomorphism classes of}\ G_{\phi_\alpha} \}\right.
$$
Consider the following cases:\\ \\
1) Let $a_1=(\alpha_1,\alpha_2,\alpha_3)$ with
$\alpha_1=[t_2,t_3],\ \alpha_2=[t_1,t_3],\ \alpha_3=[t_1,t_2]$.
Then \begin{align*} & P_{a_1}(y_1,\dots,
y_6)=y_1^2y_4-y_2^2y_5+y_3^2y_6-y_1y_2y_3,\\
& Hes(P_{a_1})(y_1,\dots, y_6)=64y_1^2y_2^2y_3^2.
\end{align*}
2) Let $a_2=(\alpha_1,\alpha_2,\alpha_3)$ with
$\alpha_1=[t_2,t_3],\ \alpha_2=[t_2,t_3],\ \alpha_3=[t_1,t_2]$.
Then
\begin{align*} &
P_{a_2}(y_1,\dots,
y_6)=y_1^2y_4-y_2y_5y_3+y_3^2y_6-y_1y_3^2,\\
& Hes(P_{a_2})=16y_1^2y_3^4.
\end{align*}
3) Let $a_3=(\alpha_1,\alpha_2,\alpha_3)$ with
$\alpha_1=[t_2,t_3],\ \alpha_2=[t_2,t_3],\ \alpha_3=[t_1,t_2]$.
Then
\begin{align*}
& P_{a_3}(y_1,\dots, y_6)
=y_1y_4y_3-y_2y_5y_3+y_3^2y_5-y_3^3,\\
& Hes(P_{a_3})=4y_3^6.
\end{align*}
Clearly, all these three cases define non-$k$-equivalent
polynomials. Hence, the corresponding groups and Lie algebras are
not isomorphic.
\begin{prop}
The Lie algebras $L(G_{a_1}), L(G_{a_2}), L(G_{a_3})$ are pairwise
non-isomorphic.
\end{prop}

\noindent{\bf Remark.} The invariant polynomials for the two-step
nilpotent quotients of Lie algebras of the groups
$Cb_4^+/Z(Cb_4^+)$ and $P_4/Z(P_4)$ define cubic forms which lie
on Hilbert's null-cone, hence they do not differ from each other
by an invariant of a Hessian type\footnote[1]{Authors thank V.L.
Popov for helping analyzing these cubic forms and for this
remark}.

 \section{On non-linearity of certain automorphism groups}
\vspace{.5cm}
\subsection{Poison group}
Let $G$ be a group. Consider the group $\mathcal H(G)$, defined as
an HNN-extension:
$$
\mathcal H(G)=\langle G\times G,t\ |\ t(g,g)t^{-1}=(g,1),\ g\in
G\rangle.
$$
It is easy to see that the group $\mathcal H(G)$ can be presented
as a semi-direct product
$$
\mathcal H(G)\simeq (G*\mathbb Z)\rtimes G,
$$
where the action of $G$ on $G*\mathbb Z=G*\langle t\rangle$ is
defined by
$$
g: g_1t^{e_1}g_2t^{e_2}\dots g_kt^{e_k}\mapsto
g_1{(tg^{-1})}^{e_1}g_2{(tg^{-1})}^{e_2}\dots
g_k{(tg^{-1})}^{e_n},\ g,g_i\in G, e_i\in\mathbb Z.
$$

For example, if $G=\mathbb Z$, we get the following group
$$
\mathcal H(G)=\langle a,b\ |\ [a,[b,a]]=1\rangle.
$$
It follows from elementary Tietze moves that
$$
H(\mathbb{Z}) = \langle a, a', t\ |\ [a, a'] = 1,~~~ t a a' t^{-1}
= a \rangle = \langle a, a', t\ |\ [a, a'] = 1,~~~ a' = a^{-1}
t^{-1} a t = [a, t] \rangle =
$$
$$
= \langle a, t\ |\ [a, [a, t]] = 1\rangle,
$$
where $a =(a, 1), a' = (1, a) \in G \times G.$

Denote by $\mathcal{NAF}$ the class of
nilpotent-by-abelian-by-finite groups. The following result due to
Formanek and Procesi presents a way of constructing  non-linear
groups.
\begin{theorem}\cite{FP}\label{fp}
If $G\notin \mathcal{NAF},$ then the group $\mathcal H(G)$ is non-
linear.
\end{theorem}
The simplest example of a group which does not lie in the class
$\mathcal{NAF}$, clearly, is a free non-cyclic group. Thus, for
$G=F_2,$ the group $\mathcal H(F_2)$, called a {\it poison group},
is non-linear. This fact plays an important role in the proof of
non-linearity of groups $Aut(F_n)$ for $\geq 3$; the poison group
is a subgroup in $Aut(F_n),\ n\geq 3.$ This statement can be
generalized.
\begin{theorem}
Let $G\notin \mathcal{NAF}$, then $Aut(G*\mathbb Z)$ is
non-linear.
\end{theorem}
\begin{proof}
We will realize the group $\mathcal H(G)$ as a subgroup in
$Aut(G*\mathbb Z)$ and the statement will follow from Theorem
\ref{fp}. Elements of the subgroup  $G\times G$ in $\mathcal H(G)$
we will denote as $(g,g'),$ i.e. we put dash for elements from the
second copy of $G$. The group $G*\mathbb Z$ we will describe in
terms of generators $g\in G$ and a free generator $t$. Consider
the homomorphism
$$
f: \mathcal H(G)\to Aut(G*\mathbb Z),
$$
given by setting
$$
f: g\mapsto i_g,\ g\in G,\ t\mapsto i_t,\ g'\mapsto s_{g'},\ g'\in
G,
$$
where $i_g$ is the conjugation by $g$, $i_t$ is the conjugation by
 $t$, $s_{g'}\ (g'\in G)$ is the automorphism of $G*\mathbb Z$
 acting trivially on $G$ and sending the element $t$ to the element $tg'^{-1}$.
It can be checked that $f$ is a group homomorphism. Every element
of the group $\mathcal H(G)$ can be written without 'dash'
elements, since $t g g' t^{-1} = g$ and, therefore, $g' = [g, t]$.
Hence, in the case of the existence of a non-trivial kernel of
$f$, there is an element of $G*\mathbb Z$, acting trivially by
conjugation, i.e. lying in the center of $G*\mathbb Z$. However,
any non-trivial free product has a trivial center, therefore,  $f$
is a monomorphism. Thus $\mathcal H(G)$ is a subgroup of
$Aut(G*\mathbb Z)$. Theorem \ref{fp} implies that the group
$Aut(G*\mathbb Z)$ is non-linear.
\end{proof}

Clearly, one can consider different embeddings of groups $\mathcal
H(G)$ in corresponding automorphism groups. Consider the case
$G=F_2$.

Let $F_3$ be a free group with basis $x_1,x_2,x_3$ and
$a_1,a_2,a_3$ some elements of $F_3$, such that the subgroup
$\langle a_1,a_2,a_3\rangle$ is free of rank 3. Define
automorphisms $\alpha_i,\ i=1,2,3$ of the group $F_3$ as
conjugation with $a_i$. The following statement can be checked
straightforwardly:
\begin{prop}\label{pre32}
Let $\phi_1,\phi_2\in Aut(F_3)$ be automorphisms which satisfy the
following conditions:
\begin{align*}
& \phi_1(a_1)=\phi_2(a_1)=a_1,\\
& \phi_1(a_2)=\phi_2(a_2)=a_2,\\
& \phi_1(a_3)=a_3a_1,\ \phi_2(a_3)=a_3a_2.
\end{align*}
Then the subgroup of $Aut(F_3)$, generated by elements
$\alpha_1,\alpha_2,\alpha_3,\phi_1,\phi_2$ is isomorphic to the
poison group $\mathcal H(F_2)$.
\end{prop}
As an example lets take $a_3=x_3$, with $a_1$ and $a_2$ arbitrary
elements of $\langle x_1,x_2\rangle,$ which do not lie in one
cyclic subgroup. Define
$$ \phi_{1} : \left\{
\begin{array}{ll}
x_{i} \mapsto x_{i}, & i = 1, 2 \\
x_{3} \mapsto x_{3}a_1,
\end{array} \right.
$$
$$
\phi_{2} : \left\{
\begin{array}{ll}
x_{i} \mapsto x_{i}, & i = 1, 2 \\
x_{3} \mapsto x_{3}a_2.
\end{array} \right.
$$
The conditions \ref{pre32} can be checked straightforwardly. Then
the subgroup of $Aut(F_3)$, generated by the elements $\alpha_i,
\phi_j,\ i=1,2,3,j=1,2$ is isomorphic to $\mathcal H(F_2)$. In
particular, in the case $a_1,a_2\in \gamma_2(F_2)$, we have the
following
\begin{theorem}\label{cor_on_lineal}
The group $IA(F_3)$ contains a subgroup isomorphic to $\mathcal
H(F_2)$, and hence, $IA(F_3)$ is not linear.
\end{theorem}

Observe also that the poison group $\mathcal H(F_2)$ is residually
finite. It follows from Baumslag's theorem which states that every
finitely generated subgroup of an automorphism group of a
residually finite group is itself residually finite. Also observe
that the non-linearity of the poison group can be used for
construction of other non-linear groups given by commutator
relations. For example, the group $\mathcal H(F_2)$ contains the
following normal subgroup of index 2:
\begin{multline*} H=\langle x_1,x_2,x_3,x_4,x_5\ |\
[x_1,x_3]=[x_2,x_4]=[x_1^{x_5},x_3]=[x_2^{x_5},x_4]=1,\\
[x_1x_3,x_2]=[x_2x_4,x_1]=[x_1^{x_5}x_3,x_4]=[x_2^{x_5}x_4,x_3]=1
\rangle,
\end{multline*}
which is of cause non-linear.

\section{Questions}
\begin{enumerate}
\item Describe the Lie algebra of the group $Cb_n,\ n\geq 2$.

\item Let $L_n$ be a free Lie algebra with $n$ generators, $n\geq
3$. Does the group $Aut(L_n)$ contain the poison group as a
subgroup?

\item Are the groups $Cb_n^+$ linear for $n\geq 3$?

\item Define the chain of subgroups
$$
{\rm Aut}(F_n) = {\rm IA}_n^1 \geq {\rm IA}_n^2 \geq {\rm IA}_n^3
\geq \ldots,
$$
where ${\rm IA}_n^k,$ $k \geq 1$ is the subgroup of ${\rm
Aut}(F_n)$, which consists of automorphisms acting trivially
modulo the $k$-th term of the lower central series of $F_n$. This
chain was introduced in \cite{A}. For which $k \geq 3,$ $n \geq 3$
the groups ${\rm IA}_n^k$ are non-linear?

\item Do the groups $Cb_n$ contain the poison group as a subgroup
for $n\geq 3$?

\end{enumerate}

 \vskip 24pt
\noindent {\it Acknowledgements.} The authors thank Vladimir
Vershinin for helpful comments and suggestions and Alexandra
Pettet for sending the preprint where the non-linearity of the
group $IA(F_5)$ is proved. 

 \end{document}